\documentclass[12pt,reqno]{amsart}
\usepackage{amssymb,delarray}
\usepackage{amsfonts,amscd}
\usepackage{epsfig}
\usepackage[all]{xy}

\def\leq{\leqslant}
\def\geq{\geqslant}

\newtheorem{thm}{Theorem}
{Lemma}
\newtheorem{prop}
{Proposition}
{Claim}
\newtheorem{df}
{Definition}
\newtheorem{cor}
{Corollary}
\newtheorem{rem}
{Remark}
{Question}
\newtheorem{ex}{Example}
\newtheorem{ex-thm}{Theorem-Example}
\newtheorem{conj}{Conjecture}
{\catcode`\@=11
\gdef\n@te#1#2{\leavevmode\vadjust{%
 {\setbox\z@\hbox to\z@{\strut#1}%
  \setbox\z@\hbox{\raise\dp\strutbox\box\z@}\ht\z@=\z@\dp\z@=\z@%
  #2\box\z@}}}
\gdef\leftnote#1{\n@te{\hss#1\quad}{}}
\gdef\rightnote#1{\n@te{\quad\kern-\leftskip#1\hss}{\moveright\hsize}}
\gdef\?{\FN@\qumark}
\gdef\qumark{\ifx\next"\DN@"##1"{\leftnote{\rm##1}}\else
 \DN@{\leftnote{\rm??}}\fi{\rm??}\next@}}

\begin{document}
\title{  Theorems of Chisini type}

\author{{ Vik.S. Kulikov}}

\address{Steklov Mathematical Institute of Russian Academy of Sciences, Moscow, Russia}
 \email{kulikov@mi-ras.ru}

\date{}

\maketitle

\def\st{{\sf st}}

\quad \qquad \qquad

\begin{abstract}
A text of the talk given by the author on  Conference on Algebra, Algebraic Geometry  and Number Theory dedicated to  the  100th anniversary of I.R. Shafarevich,  Moscow,  June 5 -- 9, 2023.
\end{abstract}

\subsection{Notations and definitions.}
Denote by $S$ a {\it smooth irreducible projective} surface defined over the field of complex numbers $\mathbb C$
and denote by
$$B:=B_f=\{ y\in \mathbb P^2\, \mid\, \sharp f^{-1}(y)<n\}$$
the {\it  branch curve} of  a finite 
map $f:S\to\mathbb P^2$, $\deg f=n$.

A divisor $R_f$ given in affine neighborhoods $U\subset S$ by
$$J(f):= \det \left(\begin{array}{cc} \frac{\partial u}{\partial z} & \frac{\partial u}{\partial w} \\ \frac{\partial v}{\partial z} & \frac{\partial v}{\partial w}\end{array}\right)= 0, $$
is called the {\it  ramification divisor } of $f$, where $(z,w)$ are local
coordinates in  neighborhoods $U\subset S$, 
$(u,v)$ are local coordinates
in $V=f(U)\subset \mathbb P^2$, and $f$ is given by $u= h_1(z,w)$,\,  $v= h_2(z,w).$

The curve ${ R=R_{f,red}}\subset S$ is called the {\it ramification curve } of $f$.
We have $f^{-1}(B)=R\cup C$, where $C$ is an additional curve to $R$  up to the proper inverse image of the branch curve.

Note that $f:S\setminus f^{-1}(B)\to \mathbb P^2\setminus B$ is unramified $n$-sheeted cover.

\subsection{Generic covers.}
\begin{thm} \label{Thm1} {\rm (Enriques, Chisini, Ciliberto, Flamini, see, for example \cite{C-F})}.  The restriction  $f:=\text{pr}_{\mid S}:S\to \mathbb P^2$ to $S$ of a linear projection
$\text{pr}: \mathbb P^N\to\mathbb P^2$ generic with respect to an embedding  $S$ in $\mathbb P^N$, has
the following properties:
\begin{itemize} \item[$(\text{i})$]  $f$ is a finite cover, $\deg f=\deg S$,
\item[$(\text{ii})$] the branch curve $B_f$ is irreducible and the points $p\in \text{Sing}\, B_f$  are only ordinary nodes and ordinary cusps,
\item[$(\text{iii})$] $f_{\mid R}: R\to B_f$ has degree $1$,
\item[$(\text{iv})$]  the ramification curve $R$  is smooth and reduced. \end{itemize} \end{thm}

A finite map  $f:S\to\mathbb P^2$ is called a {\it generic cover} if $f$ satisfies conditions $(\text{i})$ -- $(\text{iv})$.

\newpage Finite covers $f_i:S_i\to\mathbb P^2$, $i=1,2$, are {\it equivalent ($f_1\sim f_2$)}
\newline \phantom{a} if there exists an isomorphism $h:S_1\to S_2$ such that  $f_1=f_2\circ h$,

\begin{picture}(200,80)
\put(135,60){$S_1$}
\put(230,60){$S_2$}
\put(150,63){\vector(1,0){75}}
\put(185,10){$\mathbb P^2$}
\put(183,67){$h$}
\put(151,32){$f_{\text{\tiny 1}}$}
\put(217,32){$f_{\text{\tiny 2}}$}
\put(148,56){\vector(1,-1){35}}
\put(230,56){\vector(-1,-1){34}}
\end{picture}

O. Chisini in \cite{Ch} (see also \cite{Cat}) proved that two generic covers with the same branch curve are equivalent "under some suitable conditions on generality".
\newline
{\bf Chisini Conjecture.} {\it  Let $B\subset \mathbb P^2$ be the branch curve of generic covers $f_1:S\to \mathbb P^2$  and $f_2:S\to \mathbb P^2$.
If $\max (\deg f_1, \deg f_2)\geq 5$, then $f_1\sim f_2$.}

\begin{ex} Let $B\subset \mathbb P^2$ be the dual curve to a smooth cubic $\hat B\subset \mathbb P^2$.
 There exist four non-equivalent generic covers $f_i:S_i\to\mathbb P^2$, $i=1,\dots, 4$, branched in $B$, $\deg f_1=3$ and $\deg f_i=4$, $i\geq 2$.\end{ex}

Let $B$ be the branch curve of a generic cover $f:S\to\mathbb P^2$, then it is easy to show that degree of $B$ is even, $\deg B=2d$. Denote by $g$  genus of $B$ and $c$  the number of cusps of $B$.

Applying Hodge index Theorem, in \cite{K1} in 1999, it was proved
\begin{thm} \label{thm2} Let $B$ be the branch curve of two  generic covers $f_1:S_1\to\mathbb P^2$  and $f_2:S_2\to\mathbb P^2$. If
$\displaystyle \max(\deg f_1,\deg f_2)> \frac{4(3d+g-1)}{2(3d+g-1)-c}$, then $f_1\sim f_2$. \end{thm}

Using Bogomolov -- Miaoka -- Yau inequality and Theorem \ref{thm2}, S. Nemirovski (\cite{N}) proved the following
\begin{thm} \label{thm3} Let $B\subset \mathbb P^2$ be the branch curve of generic covers $f_1:S\to \mathbb P^2$  and $f_2:S\to \mathbb P^2$.
If $\max (\deg f_1, \deg f_2)\geq 12$\, then\, $f_1\sim f_2$. \end{thm}

Using Theorems \ref{thm2} and \ref{thm3}, in \cite{K2}, it was proved the following
\begin{thm}  \label{thm4}  Let $B\subset \mathbb P^2$ be the branch curve of  generic  projection $f_1:S\to \mathbb P^2$  and  generic cover
$f_2:S\to \mathbb P^2$. If $\deg f_1\geq 5$ then $f_1\sim f_2$. \end{thm}

\subsection{Monodromy of finite covers.} An unramified $n$-sheeted cover
$$f:X=S\setminus f^{-1}(B)\to Y=\mathbb P^2\setminus B$$
defines a {\it  monodromy homomorphism $f_*:\pi_1(\mathbb P^2\setminus B,q)\to \mathbb S_n$}, where $\mathbb S_n$ is the symmetric group acting on the fibre $f^{-1}(q)$.
The group $G_f:= \text{im}\, f_*\subset \mathbb S_n$ is called the {\it  monodromy group} of  $f$.

Conversely,  if a homomorphism $f_*:\pi_1(\mathbb P^2\setminus B,q)\to \mathbb S_n$ is given, then it defines an unramified $n$-sheeted cover
 $f:S\setminus f^{-1}(B)\to \mathbb P^2\setminus B$.
\newline \\
{\bf Theorem} {\rm (Riemann - Stein, \cite{St})}.
{\it For an unramified finite cover
$\widetilde f: X\to Y=\mathbb P^2\setminus B,$ there exist an uniquely defined  extension  $X\subset S$, where $S$ is a normal surface, and a finite holomorphic map
$f:S\to \mathbb P^2$ such that $f_{|X}=\widetilde f$. The surface $S$ is irreducible if and only if the monodromy group  $G_f$ is a transitive subgroup of $\mathbb S_n$.} \\

Denote by $V_p\subset \mathbb P^2$ a sufficiently small  neighborhood of a point $p\in B\subset \mathbb P^2$
bi-holomorphic to the ball $\mathbb B_{\varepsilon}=\{ (u,v)\in\mathbb C^2\mid \sqrt{u\overline u+v\overline v}<\varepsilon^2\}$, $p=(0,0)$.
If $\varepsilon\ll 1$, then $(B,p):=V_p\cap B$ is called a { \it germ} of the curve $B$ at the point $p$. \\

\noindent {\bf Theorem--Definition.} {\it Let  $V_p\simeq \mathbb B_{\varepsilon}$ be as above. Then the  group
$$\pi_1^{loc} (B,p):=\pi_1 (\mathbb B_{\varepsilon}\setminus B)$$
does not depend on $\varepsilon$ if $\varepsilon\ll 1$. The group $\pi_1^{loc} (B,p)$  is called the {\it local  fundamental group} of the curve germ $(B,p)$.}

\subsection{The local monodromy.} Let $B=B_f$ be  the branch curve of a finite cover $f:S\to\mathbb P^2$, $\deg f=n$, and let
 $U_p:=f^{-1}(V_p)\subset S$. The restriction  $f_p=f_{\mid U_p}:U_p\to V_p$ of the cover $f$ to $U_p$ is called a {\it germ} of the cover $f:S\to\mathbb P^2$ at the point $p$.

The embedding $\iota:V_p\hookrightarrow \mathbb P^2$ defines a homomorphism
$$\iota_*:\pi_1^{loc}(B,p)=\pi_1(V_p\setminus B,q)\to \pi_1(\mathbb P^2\setminus B,q)$$
and the {\it local monodromy homomorphism}
$$f_{p*}:=f_*\circ \iota_*:\pi_1^{loc}(B,p)\to\mathbb S_n.$$
The group { $G_{f,p}$}$:=\text{im}\, f_{p*}$ is called the {\it  local monodromy group} of $f:S\to\mathbb P^2$ at $p$.

The inverse image  $U_p:=f^{-1}(V_p)= \bigsqcup_{j=1}^{m_p}${ $U_{p,j}$}$\subset S$ is a disjoint union of connected surfaces $U_{p,j}$, $U_{p,j_1}\cap U_{p,j_2}=\emptyset$ for $j_1\neq j_2$.
The restriction $f_{p,j}:=f_{\mid U_{p,j}}:U_{p,j}\to V_p$, $\deg f_{p,j}=n_j$, of $f_p$ to $U_{p,j}$ is called a {\it  connected component} of the  {\it  cover germ} $f_{p}:U_p\to V_p$ of $f:S\to\mathbb P^2$,  $n=\sum_{j=1}^{m_p} n_j$.

 Let  $(B_j,p)$ be the branch curve germ of $f_{p,j}:U_{p,j}\to V_p$, $\deg f_{p,j}>1$. We have
$(B_j,p)\subset (B,p)$ and $(B_j,p)$ is called a {\it  $f$-component} of $(B,p)$.
If $\deg f_{p,j}=n_j=1$ for some $j$, then we  say that $f_{p,j}$ is {\it branched}
at $p$ and  $B_j:=p$ is called also a {\it  $f$-component} of $(B,p)$ and, by definition, $\pi_1^{loc}(B_j,p)=\text{\bf 1}$.

For the base point $q$ of the fundamental group $\pi_1^{loc}(B,p)=\pi_1(V_p\setminus B,q)$, let us consider its inverse image
$$f^{-1}(q)=\{ q_1,\dots, q_n\}\subset U_p=f^{-1}(V_p)= \bigsqcup_{j=1}^{m_p} U_{p,j}\subset S.$$
Denote by $I_{p,j}:= f^{-1}(q)\cap U_{p,j}$,\,\, { $f^{-1}(q)=\bigsqcup_{j=1}^{m_p}I_{p,j}$ and by
$\mathbb S_{n_j}:=\mathbb S(I_{p,j})\subset \mathbb S_n$ the symmetric group acting on $I_{p,j}$ and considered as a subgroup of $\mathbb S_n$,
$$\psi_I: \mathbb S_{n_1}\times \dots \times \mathbb S_{n_{m_p}}\, \hookrightarrow\, \mathbb S_n.$$

 The embedding $\iota_j: V_p\setminus B\hookrightarrow V_p\setminus B_j$ defines an epimorphism
$$\iota_{j_*}: \pi_1^{loc}(B,p)=\pi_1(V_p\setminus B)\twoheadrightarrow \pi_1(V_p\setminus B_j)=\pi_1^{loc}(B_j,p).$$

Let  $f_{p,j*}:\pi_1^{loc}(B_j,p)\to \mathbb S_{n_j}$ be the monodromy homomorphism of $f_{p,j}:U_{p,j}\to V_p$.
Consider  homomorphisms $\varphi_{p,j}$}$=f_{p,j*}\circ \iota_{j_*}: \pi_1^{loc}(B,p)\to \mathbb S_{n_j}$ and
$$\varphi_p:=\varphi_{p,1}\times\dots\times \varphi_{p,m_p}: \pi_1^{loc}(B,p)\longrightarrow \mathbb S_{n_1}\times\dots\times\mathbb S_{n_{m_p}}.$$
Then the local monodromy homomorphism $f_{p*}$ is
$$f_{p*}=\psi_I\circ \varphi_p : \pi_1^{loc}(B,p)\to\mathbb S_n.$$

Remind that singularities $(B_1,p)$ and $(B_2,p)$ are called {\it  deformation equivalent} if
 their weighted dual graphs of the minimal resolutions of $(B_1,p)$ and $(B_2,p)$ up to divisors with normal crossings, are isomorphic.

It is well-known (\cite{W}, \cite{Z1}) that if $(B_1,p)$ and $(B_2,p)$ are  deformation equivalent, then their local fundamental groups are naturally isomorphic, $\pi_1^{loc}(B_1,p)\cong \pi_1^{loc}(B_2,p)$ and if $f_{1*}:\pi_1^{loc}(B_1,p)\to \mathbb S_{n_1}$ is the monodromy homomorphism of
a finite cover $f_{1}: U_1\to V_p$, where $U_1$ is a smooth surface, then (\cite{K?}) the monodromy homomorphism $f_{2*}:\pi_1^{loc}(B_2,p)\cong \pi_1^{loc}(B_1,p)\stackrel{f_{1*}}{\longrightarrow} \mathbb S_{n_1}$ defines a finite cover $f_{2}: U_2\to V_p$ of $V_p$ in which $U_2$ is a smooth surface.

By van Kampen -- Zariski Theorem, the fundamental group $\pi_1(\mathbb P^2\setminus B,q)$ is generated by $\deg B$, so called,  {\it geometric
generators,} i.e elements $\gamma\in\pi_1(\mathbb P^2\setminus B,q)$ represented by simple loops around the curve $B$.
The local fundamental group $\pi_1^{loc} (B,p)$ is also generated by $\mu_p(B)$ geometric generators $\gamma\in \pi_1(V_p\setminus B,q)$, where $\mu_p(B)$ is
the multiplicity of the curve $B$ at the point $p\in B$.
\subsection{The case when $B$ is an irreducible curve} Next, we assume that $B$ is an {\it irreducible} curve.
In this case the geometric generators $\gamma\in \pi_1(\mathbb P^2\setminus B,q)$ are conjugated in $\pi_1(\mathbb P^2\setminus B,q)$ and
if $f_*:\pi_1(\mathbb P^2\setminus B,q)\to \mathbb S_{n}$ is the monodromy homomorphism of a finite cover $f: S\to \mathbb P^2$, then $f_*(\gamma)\in \mathbb S_n$ are conjugated in $\mathbb S_n$ for geometric generators $\gamma\in \pi_1(\mathbb P^2\setminus B,q)$.

Denote by $\overline{r}_f=(r_1,\dots ,r_k)$  the {\it  cyclical type} of permutation $f_*(\gamma)$, i.e., the collection of lengths of nontrivial cycles included in the
 factorization of $f^*(\gamma)$ into the product of disjoint cycles.
We call $\overline r_f$ the {\it  monodromy cyclical type} of the cover $f:S\to\mathbb P^2$,
the collection $\mathcal M_f=\{ f_{p*}:\pi_1^{loc}(B_f,p)\to\mathbb S_n \mid p\in \text{Sing}\, B_f\}$ is called
the {\it  local monodromy data} and a triple   $\text{pas}(f)=(B_f,\overline r_f,\mathcal M_f)$  is called the {\it  passport} of the cover $f:S\to \mathbb P^2$.

\subsection{Definition of the concept of a Chisini Theorem.} The next part of my talk is an overview of results related to answers to the following

\noindent {\bf Question.} {\it  For which finite covers $f:S\to\mathbb P^2$ do their passports uniquely determine these covers up to equivalence}? \\

Denote by $\mathcal F_{\overline r, \mathcal T}$ a set of finite covers $f:S\to\mathbb P^2$ such that $\overline r_f=\overline r$ and the singularity types of the $f$-components $(B_p,j)$ of $(B_f,p)$, $p\in \text{Sing}\, B_f$, belong to a set of singularity types $\mathcal T$ of plane curve germs.
Let $\mathcal F_{\overline r}=\displaystyle \bigcup_{\mathcal T}\mathcal F_{\overline r, \mathcal T}$.

A statement  is called {\it Chisini Theorem  with constant $\frak n$ for the covers belonging to a set $\mathcal F_{\overline r,\mathcal T}$} if it claims that  {\it there exists a constant $\frak n=\frak n(\mathcal F_{\overline r,\mathcal T})\in\mathbb N$ such that if
 $$\text{pas}(f_1)=\text{pas}(f_2)\,\,\, \text{and}\,\, \max(\deg f_1,\deg f_2)\geq \frak n$$
 for covers $f_1$ and $f_2\in\mathcal F_{\overline r, \mathcal T}$, then $f_1\sim f_2$.} \\

The following Theorem is well-known.
\begin{thm} \label{thm0} For all $r\geq 2$, Chisini Theorem with constant $\frak n=1$ holds for finite covers belonging to the set $\mathcal F_{(r),\emptyset}$. \end{thm}

Let  $\mathcal Gen =\mathcal F_{(2),\{ A_0,A_2\}}$ be the set of generic covers and $\mathcal Gen_{\text{pr}}\subset\mathcal F_{(2),\{ A_0,A_2\}}$ the set of generic projections.

Theorems 2 and 3 are examples of Chisini Theorems with constant  $\frak n=12$ for the covers belonging to the set $\mathcal Gen$ and\, with constant $\frak n=5$ for the covers
belonging to the set $\mathcal Gen_{\text{pr}}$.

\begin{prop} \label{prop1} For each $m>1$, there is no a constant $\frak n\in\mathbb N$ for  which Chisini Theorem with constant $\frak n$ would be fulfilled
 for the covers belonging to $\mathcal F_{(\underbrace {2,\dots,2}_{m}), \{ A_0,A_2\}}$.\end{prop}
\proof Let $S_{o}$ be an  Abelian surface whose endomorphism ring is $\mathbb{Z}$,
then, for example, for each prime number $m$, there are $\frac{m^{4}-1}{m-1}$ finite et$\hat{\text{a}}$le
cyclic covers $f_i:S_i\to S_{o}$ of  degree $m$ such that $S_{i_1}$ and $S_{i_2}$ are not isomorphic
for $i_1\neq i_2$. Therefore if $B$ is the branch curve of a generic projection
$f_o:S_{o}\to \mathbb P^2$, then $B$ is the branch curve of at least $\frac{m^{4}-1}{m-1}$  covers
$\{ \widetilde f_i=f_o\circ f_i:S_i\to\mathbb P^2\} \in \mathcal F_{(\underbrace {2,\dots,2}_{m}), \{ A_0,A_2\}}$ of degree $m\cdot\deg f_o$.

Note that for any $N\gg 1$ we can take a generic projection $f_o:S_o\to\mathbb P^2$ of $\deg f_o>N$. \qed \\

Proposition \ref{prop1} shows that there is no hope that Chisini Theorems hold for monodromy cyclical types that do not consist of a single cycle. Note also that by now almost nothing is known about theorems of Chisini type in the case when the monodromy cyclical type is $\overline r_f =(r)$ with $r\geq 3$, except, perhaps,  Theorem \ref{thm0}. Therefore, further we will limit ourselves to considering the case when $r=2$.

\subsection{ The case of monodromy cyclical type $\overline r_f=(2)$.}
In order to be able to formulate statements related to the answers to the question formulated above, it is necessary to recall a few more definitions related to the invariants of the singularities of plane curves.

Let $(B,p)=V_p\cap B$ be a germ of a curve $B\subset \mathbb P^2$ at a point $p$.  Denote by $\delta_p(B)$ the {\it  $\delta$-invariant} of the singularity $(B,p)$,
and denote by $\mu_p(B)$ the {\it  multiplicity} of the singularity $(B,p)$ at $p$.

Let $(B,p)=(B_1,p)\cup\dots\cup (B_k,p)$ be the union of $k$ irreducible germs $(B_i,p)$ of $(B,p)$,
The number
$$c_{v,p}(B):=\sum_{i=1}^k (\mu_p(B_i)-1)$$ is called the {\it  number of virtual cusps} of the curve germ $(B,p)$ and
 $$c_{v}(B):=\sum_{p\in\text{Sing}\, B} c_{v,p}$$
 is called the {\it  number of virtual cusps} of the curve $B\subset \mathbb P^2$.
The number
$$n_{v}(B):=\delta(B)-c_{v}(B)$$
is called the {\it  number of virtual nodes}  of the curve  $B$.

Note that the numbers of virtual cusps and nodes of a plane curve and its dual curve play the same role as the numbers of ordinary casps and ordinary nods in the case of cuspidal curves, that is, if we substitute the numbers of virtual cusps and nods of a curve and its dual curve in the classical Pl\"ucker formulas instead of the numbers of ordinary cusps and nods, then we get again equality in these formulas (see \cite{Ku}).

As it was mentioned above, if $f_*:\pi_1(\mathbb P^2\setminus B,q)\to \mathbb S_{n}$ is the monodromy homomorphism of $f: S\to \mathbb P^2$, then
the images $f_*(\gamma)\in \mathbb S_n$ are transpositions in $\mathbb S_n$ for geometric generators $\gamma\in \pi_1(\mathbb P^2\setminus B,q)$.
Consequently,  $G_f=\mathbb S_n$.

Let $(B,p)=V_p\cap B$ be a germ of $B$ at $p\in \text{Sing}\, B$, $f_{p,j}:U_{p,j}\to V_p$, $\deg f_{p,j}=n_j$,  a connected component  of the cover germ $f_{p}:U_p\to V_p$ of $f:S\to\mathbb P^2$, and $(B_j,p)\subset (B,p)$ is the $f$-component of $(B,p)$.
Since fundamental group $\pi_1^{loc}(B_j,p)$ is generated by $\mu_p(B_j)$ geometric generators, then $G_{f,p}=\mathbb S_{n_j}$, $n_j\leq \mu_p(B_j)+1$.
 A connected component $f_{p,j}:U_{p,j}\to V_p$ is called {\it non-degenerate} if $n_j=\mu_p(B_j)+1$.

\subsection{Extra property.} Consider two copies of a connected component $f_{p,j}:U_{p,j}\to V_p$ of a cover germ $f_{p}:U_p\to V_p$ of $f:S\to\mathbb P^2$. Denote by
$$f_i:=f_{p,j} :U_i:=U_{p,j}\to V_p,$$ $\deg f_i=n_j$, $f_i^{-1}(B_j)=R_i\cup C_i\subset U_i$, $i=1,2$, and
consider  the fibre product
$$U_1\times _{V_p}U_2=\{ \, (x,y)\in U_1\times U_2\, \,  \mid \, \, f_1(x)=f_2(y)\, \, \}. $$
 Let { $\widetilde X$}$=\widetilde{U_1\times _{V_p}U_2}$ be the normalization of $U_1\times _{V_p}U_2$.
 Denote by
 $$g_{\text{\tiny 1}}:\widetilde X\to U_1,\qquad  g_{\text{\tiny 2}}:\widetilde X\to U_2,\qquad  g_{\text{\tiny 1,2}}:\widetilde X\to V_p$$
the corresponding natural morphisms,  $\deg g_{\text{\tiny 1}}=\deg g_{\text{\tiny 2}}=n_j$,   $\deg g_{\text{\tiny 1,2}}=n_j^2$.

Note that $\widetilde X=\widetilde X'\bigsqcup \widetilde X''$ is a disjoint union of two connected components $\widetilde X'$ and $\widetilde X''$,
$\deg g_{i\mid \widetilde X'}=n_j-1$ and { $\deg g_{i\mid \widetilde X''}=1$.

\newpage We have the following commutative diagram.

\begin{picture}(350,180)
\put(112,90){$p_ {\text{\tiny 1}}\in U_{\text{\tiny 1}}$}
\put(230,90){$U_{\text{\tiny 2}}\ni p_{\text{\tiny 2}}$}
\put(184,42){$V_p\ni p$}
\put(183,146){$\widetilde X'\ni p_ {\text{\tiny 1,2}}$}
\put(180,142){\vector(-1,-1){38}}
\put(197,142){\vector(1,-1){38}}
\put(152,129){$ g_{\text{\tiny 1}}$}
\put(215,129){$ g_{\text{\tiny 2}}$}
\put(150,63){$f_{\text{\tiny 1}}$}
\put(217,63){$f_{\text{\tiny 2}}$}
\put(148,85){\vector(1,-1){32}}
\put(230,85){\vector(-1,-1){32}}
\put(174,20){$\text{Fig. 1}$}
\end{picture}

Note that $p_{\text{\tiny 1,2}}=g_1^{-1}(f_1^{-1}(p))\cap \widetilde X'$ is, possibly,  the only singular point of $\widetilde X'$.

Let us add the inverse images of the curve germ $(B_j,p)$ into the commutative diagram depicted in Fig 1. We obtain a commutative diagram depicted in Fig 2 in which
$$\widetilde R\cup \widetilde C_{\text \tiny 1}=g_{\text \tiny 1}^{-1}(R_{\text \tiny 1}),\qquad
\widetilde R\cup \widetilde C_{\text \tiny 2}=g_{\text \tiny 2}^{-1}(R_{\text \tiny 2}),$$
and
$$\widetilde R=g_{\text{\tiny 1}}^{-1}(R_{\text{\tiny 1}})\cap g_{\text{\tiny 2}}^{-1}(R_{\text{\tiny 2}}).$$

\begin{picture}(350,200)
\put(160,90){$U_{\text{\tiny 1}}$}
\put(80,92){$p_{\text{\tiny 1}}\in C_{\text{\tiny 1}}\cup R_{\text{\tiny 1}}\subset$}
\put(250,90){$U_{\text{\tiny 2}}$}
\put(265,92){$\supset R_{\text{\tiny 2}}\cup C_{\text{\tiny 2}} \supset p_{\text{\tiny 2}}$}
\put(202,42){$V_p\ni p$}
\put(88,42){$V_p\cap B_f\supset (B_j,p)\subset $}
\put(204,150){$\widetilde X'$}
\put(115,150){$p_{\text{\tiny 1,2}}\in\widetilde R \cup \widetilde C_{\text{\tiny 1}}\subset$}
\put(225,150){$\supset \widetilde R \cup \widetilde C_{\text{\tiny 2}}\ni p_ {\text{\tiny 1,2}}$}
\put(200,142){\vector(-1,-1){35}}
\put(160,145){\vector(-1,-2){20}}
\put(257,145){\vector(1,-2){20}}
\put(217,142){\vector(1,-1){35}}
\put(125,88){\vector(1,-1){33}}
\put(175,129){$g_{\text{\tiny 1}}$}
\put(234,129){$g_{\text{\tiny 2}}$}
\put(169,65){$f_{\text{\tiny 1}}$}
\put(238,65){$f_{\text{\tiny 2}}$}
\put(168,85){\vector(1,-1){30}}
\put(250,85){\vector(-1,-1){30}}
\put(197,20){$\text{Fig. 2}$}
\end{picture}

Denote by $\nu : X'\to\widetilde X'$ a resolution of singularity of $\widetilde X'$. We obtain the following commutative diagram
in which
$\overline R=\nu^{-1}(\widetilde R)$ and $\overline C_1=\nu^{-1}(\widetilde C_1)$ are the proper inverse images of curve germs $\widetilde R$ and $\widetilde C_1$ and  $(g_1\circ\nu)^*(R_1)=\overline R+\overline C_1+E$ is the inverse image of divisor $R_1\in \text{Div}\, U_1$,
$\text{Supp}(E)\subset \nu^{-1}(p_{\text{\tiny 1,2}})$.

\begin{picture}(350,250)
\put(180,82){$U_{\text{\tiny 1}}$}
\put(125,83){$p_{\text{\tiny 1}}\in R_{\text{\tiny 1}}\subset$}
\put(280,82){$U_{\text{\tiny 2}}$}
\put(229,32){$V_p\ni p$}
\put(120,32){$V_p\cap B_f\supset (B_j,p)\subset $}
\put(232,140){$\widetilde X'\ni p_ {\text{\tiny 1,2}}$}
\put(166,141){$\widetilde R \cup \widetilde C_{\text{\tiny 1}}\,\, \subset$}
\put(230,137){\vector(-1,-1){42}}
\put(180,136){\vector(-1,-2){20}}
\put(240,137){\vector(1,-1){42}}
\put(198,119){$g_{\text{\tiny 1}}$}
\put(233,220){$X'\ni \nu^{-1}( p_ {\text{\tiny 1,2}})$}
\put(237,215){\vector(0,-1){60}}
\put(240,185){$\nu$}
\put(262,119){$g_{\text{\tiny 2}}$}
\put(155,220){$\overline R + \overline C_{\text{\tiny 1}}+ E\subset$}
\put(196,55){$f_{\text{\tiny 1}}$}
\put(182,216){\vector(0,-1){60}}
\put(265,55){$f_{\text{\tiny 2}}$}
\put(190,78){\vector(1,-1){33}}
\put(280,80){\vector(-1,-1){35}}
\put(200,5){$\text{Fig. 3}$}
\end{picture}

\begin{df} {\rm (\cite{Kul})} A germ  $f_{\text{\tiny 1}}=f_{p,j} : U_1=U_{p,j}\to V_p$ has {\it extra property} over the point $p\in\text{Sing}\, B$
 if either $\deg f_{\text{\rm \tiny 1}}\leq 2$, or if $\deg f_{\rm \text{\tiny 1}}> 2$, then there is a presentation of $E$ in the form $E=E_{\overline R}+E_{\overline C_1}$ such that
 $$(\overline R+E_{\overline R},\overline C_1+E_{\overline C_1})_{X'}\leq  2 \delta_{p_{\rm \text{\tiny 1}}}(R_1)+c_{v,p}(B_j).$$
\end{df}

Denote by  $\mathcal T_{Ch_{\text{\tiny 12}}}$ the set of {\it singularity types} of $f$-components 
$(B_j,p)$ of {\it connected non-degenerate} cover germs
$f_{p,j}:U_{p,j}\to V_p$ having { extra property}.
\begin{thm} {\rm (\cite{K7}, \cite{Kul})}  Chisini Theorem with constant $\frak n=12$ holds for covers $f:S\to \mathbb P^2$ belonging to $\mathcal F_{(2),\mathcal T_{Ch_{\text{\rm \tiny 12}}}}$, i.e., if
$pas(f_1)=pas(f_2)$  and $\max(\deg f_1,\deg f_2)\geq 12$  for covers $f_1:S_1\to\mathbb P^2$ and $f_2:S_2\to\mathbb P^2$ belonging to
$\mathcal F_{(2),\mathcal T_{Ch_{\text{\rm \tiny 12}}}}$, then $f_1\sim f_2$.
\end{thm}

\begin{thm} {\rm(1)} The set of $f$-components of branch curve germs of $\text{ADE}$-singularity
\newline \phantom{Theorem 6.aaaaa} types is $\{ A_0, E_6, A_{3n+2}, n\geq 0\}$.
\newline \phantom{Theorem 6.aa} {\rm (2)}   $\{ A_0, E_6, A_{3n+2}, n\geq 0\}\subset \mathcal T_{Ch_{\text{\rm \tiny 12}}}$.
\end{thm}

\begin{cor} Chisini Theorem with constant $\frak n=12$ holds for finite  covers $f:S\to\mathbb P^2$ belonging to
$\mathcal F_{(2),\{ A_0, E_6, A_{3n+2},\,\, n\geq 0\}}$.
\end{cor}

\subsection{Dualizing covers of the plane associated with curves immersed in the plane.} \label{ex1}
Let $\iota: C\to {\mathbb P}^2$ be a morphism of smooth irreducible reduced projective curve $C$ to the projective plane such that $\iota_{\mid C} :C\to \iota(C)$ is a bi-rational morphism, $\deg \iota(C)=n\geq 2$ and let $B\subset \hat{\mathbb P}^2$ be the dual curve to the curve $\iota(C)$ (the curve $B$ consists of lines $l_p\in \hat{\mathbb P}^2$, $p\in C$, tangent to the curve $\iota(C)$). The correspondence graph between curves $\iota(C)$ and $B$ is a curve  $\check C$ in ${\mathbb P}^2\times \hat{\mathbb P}^2$ (the so-called Nash curve blowup of $\iota(C)$), which lies in the variety of incident $I=\{ (P,l)\in {\mathbb P}^2\times \hat{\mathbb P}^2\mid P\in l\}$,
$$\check C=\{ (\nu(p),l_p)\in I\mid p\in  C\, \, \text{and}\, \, l_p  \, \text{is a tangent line to}\, \, \iota(C)\,\, \text{at}\,\, \iota(p)\in C\}.$$

Let $\text{pr}_1:\mathbb P^2\times \hat{\mathbb P}^2\to \mathbb P^2$ and $\text{pr}_2:\mathbb P^2\times \hat{\mathbb P}^2\to \hat{\mathbb P}^2$ be the projections to the factors, $X'=\text{pr}_1^{-1}(\iota(C))\cap I$, and $h:X'\to \hat{\mathbb P}^2$ the restriction of the projection $\text{pr}_2$ to $X'$. Obviously, $h^{-1}(l)$ consists of the points $(P,l)\in\mathbb P^2\times \hat{\mathbb P}^2$ such that $P\in \iota(C)\cap l$ and hence, $\deg h=\deg \iota(C)=n$.

In \cite{K-d}, the notion of the {\it dualizing cover of the plane associated with the curve} $\iota(C)\subset \mathbb P^2$ was introduced. Namely, let $\nu:X\to X'$ be the normalization of the surface $X'$, then the morphism $f_{\iota(C)}=h\circ \nu:X\to \hat{\mathbb P}^2$ is called the {\it dualizing cover} of the plane associated with the curve $\iota(C)\subset \mathbb P^2$.  We have $\deg f_{\iota(C)}=n$. It is easy to see that $X$ is isomorphic to the fibre product $C\times_{\iota(C)}X'$ of the morphism $\iota: C\to \iota(C)$ and the projection $\text{pr}_1:X'\to \iota(C)$. The projection $\text{pr}_1:X'\to \iota(C)$ defines a structure of ruled surface on $X'$ and this structure induces a ruled structure on $X$ over $C$. Note that  $\widetilde {C}=\nu^{-1}(\check C)\subset X$ is a section of this ruled structure, and the image $f_{\iota(C)}(F_p)$ of a fibre $F_p$ is a line $L_{\iota(p)}\subset \hat{\mathbb P}^2$ dual to the point $\iota(p)\in \iota(C)\subset \mathbb P^2$.

\begin{thm} \label{main5} Let $\iota: C\hookrightarrow \hat{\mathbb P}^2$ be an immersion of irreducible projective curve $C$ of genus $g$ and $f_{\iota(C)}:X_{\iota(C)}\to \mathbb P^2$ the dualizing cover associated with $\iota(C)$, and let a cover $f:X\to \mathbb P^2$ belong to the set $\mathcal F_{(2)}$. If $pas(f)=pas(f_{\iota(C)})$, then the covers $f$ and $f_{\iota(C)}$ are equivalent if either $\deg \iota(C)\geq 8$ and $g\geq 1$, or $\deg \iota(C)\geq 12$ and $g=0$.
\end{thm}

Consider covers given by functions $u=z$, $v=w^{n}-nwz$. They are branched in curves given by equation $v^{n-1}-(1-n)^{n-1}u^{n}=0$.
 Denote by  $T_{n,1}$, $n\geq 2$, the singularity type  at the origin of function $v^{n-1}-(1-n)^{n-1}u^{n}$.
\begin{rem} We have $T_{2,1}=A_0$,\,\, $T_{3,1}=A_2$, and $T_{4,1}=E_6$.\end{rem}

\begin{thm} {\rm (\cite{K7}, \cite{Kul})} $\{ T_{n,1},\,\, n\geq 2\}\subset \mathcal T_{Ch_{\text{\rm \tiny 12}}}$, i.e., Chisini Theorem with constant $\frak n=12$ holds for finite covers $f:S\to\mathbb P^2$ belonging to $\mathcal F_{(2),\{ T_{n,1},\,\, n\geq 2\}}$.\end{thm}

Denote by { $T_{n,m}$} the singularity type of the branch curve germ  $B_{n,m}$ of a cover germ  $f_{p,n,m,}: U_{p,n,m}\to V_p$ given by
functions $u=z$, $v=w^n-nwz^m$.
\begin{conj} $\{ T_{n,m},\,\, n\geq 2,m\geq 1\}\subset \mathcal T_{Ch_{\text{\rm \tiny 12}}}$, i.e., Chisini Theorem with constant $\frak n=12$ holds for finite covers $f:S\to\mathbb P^2$ belonging to $\mathcal F_{(2),\{ T_{n,m},\,\, n\geq 2, m\geq 1\}}$. \end{conj}

\noindent {\bf Question.} {\it Is it true that $ 
\mathcal F_{(2)} 
=\mathcal F_{(2),\mathcal T_{Ch_{\text{\tiny 12}}}}$}?

\ifx\undefined\bysame
\newcommand{\bysame}{\leavevmode\hbox to3em{\hrulefill}\,}
\fi


\end{document}